\documentclass[11pt]{article}
\newtheorem{thm}{Theorem}[section]
\newtheorem{prop}{Proposition}[section]
\newfont{\inte}{msbm10}  
\newcommand{\Prob}{{\rm P}}
\newcommand{\EE}{{\rm E}}
\newcommand{\Var}{{\rm Var}}
\newcommand{\Stati}{{\cal S}}
\newcommand{\ds}{\displaystyle}
\newcommand{\qed}{\hfill $\diamond$}  

\newenvironment{proof}{
\begin{trivlist}
\item[\hspace{\labelsep}{\bf\noindent Proof. }]
}{\end{trivlist}}
\begin{document}
\setlength{\baselineskip}{20pt}   
\title{\huge\bf On the first-visit-time problem \\
for birth and death processes \\
with catastrophes\footnote{short title: 
``Birth-death processes in the presence of catastrophes''}
}
\author{
{\sc A.\ Di Crescenzo}$^{(1)}$, 
{\sc V.\ Giorno}$^{(1)}$,  
{\sc A.G.\ Nobile}$^{(1)}$, 
{\sc L.M.\ Ricciardi}$^{(2)}$\footnote{Corresponding author}  \\
\hfill \\
${(1)}$ \ Dipartimento di Matematica e Informatica \\ 
Universit\`a di Salerno\\ 
Via S.\ Allende, I-84081 Baronissi (SA), Italy. \\
{\normalsize Email address: \{adicrescenzo,giorno,nobile\}@unisa.it}  \\
\hfill \\
${(2)}$ \ Dipartimento di Matematica e Applicazioni \\
Universit\`a di Napoli Federico II\\ 
Via Cintia, I-80126 Napoli, Italy. \\
{\normalsize Email address: luigi.ricciardi@unina.it}
}
\date{\today}
\maketitle
\begin{abstract}
For a birth-death process subject to catastrophes, defined on the state-space 
$\Stati=\{r,r+1,r+2,\ldots\}$, with $r$ a positive integer or zero, the 
first-visit time to a state $k\in\Stati$ is considered and the Laplace transform 
of its probability density function is determined, use of which is then made 
to obtain mean and variance. The Laplace transform of the probability density function 
of the first effective catastrophe occurrence time and its expected value are also 
obtained. Some extensions to time-non-homogeneous processes are then provided. 
Finally, certain additional results concerning the determination of the steady-state 
distribution and the representation of the transition probabilities are worked out, 
while some applications to particular birth-death processes are shown in the Appendix. 


\end{abstract}
%
\section{Introduction}
Great attention has been paid in the literature to the description of the 
evolution of systems modeled via discrete state-space random processes 
such as populations evolving in random environments or queueing and service 
systems under various operating protocols. However, with a few exceptions that
will be indicated below, in these studies the systems under considerations 
are not subject to catastrophes, neither a target of investigation has 
been the time distribution when a preassigned state is reached. 
Among the relevant contributions to the area in which the present 
paper belongs, the following should be recalled: {\em (i)} The results concerning the 
distribution of the extinction time for a linear birth-death process subject to 
catastrophes and, in particular, the determination of necessary and sufficient 
conditions for population extinction to occur \cite{Br85}, \cite{Br86}; {\em (ii)} the 
evaluation of the stationary probabilities of a simple immigration-birth-death 
process influenced by total catastrophes \cite{Ky94}; {\em (iii)} the studies on the 
transient and equilibrium behaviors of immigration-birth-death process with 
catastrophes \cite{ChZh2003}, \cite{Sw97}; {\em (iv)} the analysis of linear birth-death 
processes under the influence of Poisson time-distributed catastrophes whose 
effect is to enhance the death probability  \cite{BaBuChPe89}, \cite{PePeChBa93}; 
{\em (v)} the determination of the transient and the limiting distributions 
of continuous-time Markov chains subject to catastrophes occurring according 
to renewal processes \cite{EcFa2003}; {\em (vi)} analysis of the effect of 
catastrophes in the case of $M/M/1$ queueing systems \cite{DiGiNoRi03} and 
when the number of initially present customers is random \cite{KrAr2000}. 
\par
The determination of the distribution of the time when the zero-state is reached 
for the first time, studied in some of the above papers, is clearly of interest 
within the context of population 
dynamics (cf.\ extinction problems) or in operations research (cf.\ the running out 
of stored supplies, the disappearance of a queue at a service station, etc.). 
However, no mention appears to have been made so far on the determination of the 
distribution of the time when first a discrete continuous-time Markov process of 
birth-death type visits a preassigned non-zero state when subject to catastrophes, 
i.e.\ to the fist-visit-time problem. This is the object of the present paper. 
\par
It should not pass unnoticed that the fist-visit-time problem bears particular 
relevance in contexts such as populations for which extinction does not occur 
only when population size vanishes, but is an automatic consequence of its size 
reduction below some critical values (think of bisexual populations or environmental 
effects by which very low densities random mating individuals automatically lead 
to sure extinction \cite{Ri77}). Moreover, it is conceivable that in some instances, 
due to environmental regulation effects when a system reaches too large a size 
some new phenomena may arise whose effect may be determinant on the further 
evolution of the system, such as the stop of the access to a service provider 
if a waiting room fills up, the block of an information processing system as 
a consequence of excessive service demands, cannibalism phenomena among starving 
animal populations, etc., all this both in the presence and in the absence of 
catastrophes. In all such cases, the problem arises of determining the ``statistics'' 
of the time when a generally non-zero preassigned state is for the first time 
reached by the system under consideration. 
\par
In the sequel, we shall preliminarly establish certain relations among transition 
probabilities in the presence and in the absence of absorbing states, and then reduce 
the study of the dynamics of the process in the presence of catastrophes to that 
of a process insensitive to catastrophes, whose description is in principle more 
easily achieved. Some additional results and illustration examples to particular 
birth-death processes will be briefly outlined in the Appendix. 
\par
More specifically, in Section 2 we shall analyze the first-visit time of a 
birth-death process with catastrophes $\{N(t);\, t\geq 0\}$ to any preassigned 
state $k\in\Stati=\{r,r+1,r+2,\ldots\}$, with $r$ a positive integer or zero, 
by considering separately the cases of first visit to state $r$ and to state 
$k\neq r$. We shall disclose various functional relations that allow to describe 
$N(t)$ in terms of the corresponding birth-death process $\{\widehat N(t); t\geq 0\}$, 
defined on the same state-space $\Stati$ and characterized by the same birth and 
death rates as $N(t)$, but for which catastrophes are absent. In particular, for 
$N(t)$ we obtain the Laplace transform of the probability density function (pdf) 
of the first-visit time to state $k$ and calculate its mean and variance. 
\par 
We point out that since $N(t)$ is skip-free to the right, all intermediate 
states are visited only when the first visit to state $k$ occurs from a 
starting state $j<k$. Note that, because of the assumed presence of 
catastrophes, this need not occur if $j>k$. 
\par
The problem of the first occurrence of catastrophe is faced in Section 3, 
where the Laplace transform of the catastrophe's first-occurrence time pdf 
and its mean are obtained. 
Our approach is analogous to that exploited in \cite{DiGiNoRi03} for a 
special birth-death process of interest in queueing theory. 
In Section 4 we shall obtain the expressions of the transition probabilities 
and of the density of the first-visit time of $N(t)$ to state $r$ in the 
special case when this process is time non-homogeneous and catastrophes 
occur with time-varying rates. 
In Section 5 we express the steady-state probabilities of $N(t)$ in terms 
of the Laplace transforms of the transient probabilities of $\widehat N(t)$, 
and provide a probabilistic interpretation of the transition probabilities 
of $N(t)$. Indeed, we show that $N(t)$ has the same distribution as the 
minimum between the random variable that describes its steady state, and an 
independent birth-death process whose rates suitably depend on those of $N(t)$. 
\par
Some examples of applications of our results are finally given in the Appendix 
for birth process, immigration-emigration process, immigration-death process, 
and immigration-birth-death process in the presence of catastrophes. 
\section{First-visit-time problem}
Let $\{N(t);\, t\geq 0\}$ be a birth-death process with catastrophes defined 
on the state-space $\Stati=\{r,r+1,r+2,\ldots\}$, with $r$ zero or a positive 
integer, such that transitions occur according to the following scheme: \\
(i) \ $n\to n+1$ with rate $\alpha_n$, for $n=r,r+1,\ldots$, \\
(ii) \ $n\to n-1$ with rate $\beta_n$, for $n=r+2,r+3,\ldots$, \\
(iii) \ $r+1\to r$ with rate $\beta_{r+1}+\xi$, \\
(iv) \ $n\to r$ with rate $\xi$, for $n=r+2,r+3,\ldots$. 
\par
Hence, births occur with rates $\alpha_n$, deaths with rates $\beta_n$, and 
catastrophes with rate $\xi$, the effect of each catastrophe being the 
instantaneous transition to the reflecting state $r$. For all $j,n\in\Stati$ 
and $t>0$ the transition probabilities 
$$
 p_{j,n}(t)=\Prob\{N(t)=n\,|\,N(0)=j\}
$$
satisfy the following system of forward equations:
\begin{eqnarray}
 && \hspace{-1cm} \frac{\rm d}{{\rm d}t} p_{j,r}(t) 
 =-(\alpha_r+\xi)\,p_{j,r}(t)+\beta_{r+1}\,p_{j,r+1}(t)+\xi, 
 \nonumber \\
 && \hspace{-1cm} \frac{\rm d}{{\rm d}t} p_{j,n}(t) 
 =-(\alpha_n+\beta_n+\xi)\,p_{j,n}(t)
 +\alpha_{n-1}\,p_{j,n-1}(t)   
 \label{equation:2} \\
 && \hspace{1cm}  
 +\beta_{n+1}\,p_{j,n+1}(t),  
 \hspace{2cm}  n=r+1,r+2,\ldots,
 \nonumber
\end{eqnarray}
with initial condition 
$$
 \lim_{t\downarrow 0}p_{j,n}(t)=\delta_{j,n}=\left\{
 \begin{array}{ll}
 1, & n=j \\
 0, & {\rm otherwise}.
 \end{array}
 \right.
$$
\par
Denote by $\{\widehat N(t);\, t\geq 0\}$ the time-homogeneous birth-death 
process obtained from $N(t)$ by removing the possibility of catastrophes, 
i.e.\ by setting $\xi=0$. Its transition probabilities 
$$
 \widehat p_{j,n}(t)
 =\Prob\{\widehat N(t)=n\,|\,\widehat N(0)=j\}, 
 \qquad j,n\in\Stati, \quad t\geq 0
$$
then satisfy the system of forward equations obtained from 
(\ref{equation:2}) by setting $\xi=0$, with initial condition 
$\ds\lim_{t\downarrow 0}\widehat p_{j,n}(t)=\delta_{j,n}$. 
\par
Hereafter, we shall restrict our attention to non-explosive processes 
$\widehat N(t)$, i.e.\ we shall assume that $\sum_{n=r}^{+\infty}\widehat p_{j,n}(t)=1$ 
for all $j\in\Stati$ and $t\geq 0$. 
\par
We note that some descriptors of $N(t)$ can be 
expressed in terms of the corresponding ones of $\widehat N(t)$. 
Indeed, making use of the forward equations for probabilities $p_{j,n}(t)$ 
and for $\widehat p_{j,n}(t)$, for all $j,n\in\Stati$ and $t>0$ we have: 
\begin{equation}
 p_{j,n}(t)=e^{-\xi t}\,\widehat p_{j,n}(t)
 +\xi \int_{0}^{t} e^{-\xi \tau}\,\widehat p_{r,n}(\tau)\,{\rm d}\tau.
 \label{equation:5}
\end{equation}
As an immediate consequence of Eq.\ (\ref{equation:5}), we can express the conditional 
moments of $N(t)$ in terms of those of $\widehat N(t)$ as follows:
\begin{equation}
 \EE\big[N(t)\,|\,N(0)=j\big]
 =e^{-\xi t}\,\EE\big[\widehat N(t)\,|\,\widehat N(0)=j\big] 
 +\xi \int_{0}^{t} e^{-\xi \tau}\,\EE\big[\widehat N(\tau)\,|\,\widehat N(0)=r\big]\,{\rm d}\tau. 
 \label{equation:27}
\end{equation}
Moreover, by setting 
$$
 \pi_{j,n}(\lambda)
 :=\int_0^{+\infty}e^{-\lambda t}\,p_{j,n}(t)\,{\rm d}t, 
 \qquad 
 \widehat\pi_{j,n}(\lambda)
 :=\int_0^{+\infty}e^{-\lambda t}\,\widehat p_{j,n}(t)\,{\rm d}t, 
 \qquad \lambda>0,
$$ 
from (\ref{equation:5}) it follows 
\begin{equation}
 \pi_{j,n}(\lambda)=\widehat \pi_{j,n}(\lambda+\xi)
 +\frac{\xi}{\lambda}\,\widehat \pi_{r,n}(\lambda+\xi), 
 \qquad \lambda>0.
 \label{equation:6}
\end{equation}
\par
Let us define the first-visit time of $N(t)$ to state $k\in\Stati$ as 
\begin{equation}
 T_{j,k}:=\inf\{t\geq 0:\, N(t)=k\}, 
 \qquad N(0)=j\in\Stati,\;\; j\neq k,
 \label{equation:70} 
\end{equation}
and denote its pdf by 
$$
 g_{j,k}(t)=\frac{\rm d}{{\rm d}t}\Prob\{T_{j,k}\leq t\}.
$$
Moreover, for the corresponding birth-death process $\widehat N(t)$ in the absence 
of catastrophes, we shall denote by $\widehat T_{j,k}$ the first-visit time, 
and by $\widehat g_{j,k}(t)$ its pdf. 
\par 
We shall now analyze the first-visit-time problem. To this purpose, 
$\gamma_{j,k}(\lambda)$ and $\widehat \gamma_{j,k}(\lambda)$ will denote the Laplace 
transforms of $g_{j,k}(t)$ and $\widehat g_{j,k}(t)$, respectively. Furthermore, 
for all $t\geq 0$ and $k,j,n\in\Stati$, with $j\neq k$ and $n\neq k$, let   
$$
 A_{j,n}^{\langle k\rangle}(t)
 :=\Prob\{N(t)=n, N(\tau)\neq k\;\; \forall\tau\in(0,t)\,|\,N(0)=j\}  
$$
be the $k$-avoiding transition probability of $N(t)$, and let 
$\widehat A_{j,n}^{\langle k\rangle}(t)$ denote the corresponding 
probability for $\widehat N(t)$. For all $t\geq 0$ these functions 
are related as follows: 
\begin{equation}
 \int_0^{t} \widehat g_{j,k}(\tau)\,{\rm d}\tau
 =\left\{
 \begin{array}{ll}
 1-\ds\sum_{n=k+1}^{+\infty}\widehat A_{j,n}^{\langle k\rangle}(t), & j>k\\ 
 \hfill\\
 1-\ds\sum_{n=r}^{k-1}\widehat A_{j,n}^{\langle k\rangle}(t), & j<k,
 \end{array}
 \right.
 \label{equation:68} 
\end{equation}
\begin{equation}
 \int_0^{t} g_{j,k}(\tau)\,{\rm d}\tau
 =\left\{
 \begin{array}{ll}
 1-\ds\sum_{n=r \atop n\neq k}^{+\infty} A_{j,n}^{\langle k\rangle}(t), & j>k\\ 
 \hfill\\
 1-\ds\sum_{n=r}^{k-1} A_{j,n}^{\langle k\rangle}(t), & j<k,
 \end{array}
 \right.
 \label{equation:69} 
\end{equation}
where $\int_{0}^{t} g_{j,k}(\tau)\,{\rm d}\tau$ is the probability 
that $N(t)$ enters state $k$ before time $t$, and  
$\sum_{n} A_{j,n}^{\langle k\rangle}(t)$ is the probability 
that $N(t)$ does not enter state $k$ up to time $t$.  
Hereafter, ${\cal A}_{j,n}^{\langle k\rangle}(\lambda)$ and 
$\widehat{\cal A}_{j,n}^{\langle k\rangle}(\lambda)$ will denote the 
Laplace transforms of $A_{j,n}^{\langle k\rangle}(t)$ and 
$\widehat A_{j,n}^{\langle k\rangle}(t)$, respectively.
\par
The first-visit-time problem will be discussed in three 
different cases: 
{\it (i)} first visit in state $r$, 
{\it (ii)} first visit in state $k>r$ starting from below, 
{\it (iii)} first visit in state $k>r$ starting from above. 
\subsection{First visit in state $r$}
Recalling (\ref{equation:70}), the $r$-avoiding transition probabilities of 
$N(t)$ and $\widehat N(t)$ for all $t\geq 0$  are related as follows:
\begin{equation}
 A_{j,n}^{\langle r\rangle}(t)
 =e^{-\xi t}\,\widehat A_{j,n}^{\langle r\rangle}(t), 
 \qquad j,n\in\{r+1,r+2,\ldots \}.
 \label{equation:19}
\end{equation}
Hence, the probability of a transition of $N(t)$ from $j$ to $n$ $(n>r)$ at time $t$ 
in the absence of visits to state $r$ equals the probability of the same transition 
for $\widehat N(t)$ times the probability $e^{-\xi t}$ that no catastrophe occurs 
in $(0,t)$. 
\par
The following proposition expresses $g_{j,r}(t)$ and $\gamma_{j,r}(\lambda)$ in 
terms of $\widehat g_{j,r}(t)$ and $\widehat \gamma_{j,r}(\lambda)$, respectively. 
\begin{prop}
For $j\in\{r+1,r+2,\ldots\}$ the following equations hold: 
\begin{equation}
 g_{j,r}(t)=e^{-\xi t}\,\widehat g_{j,r}(t)
 +\xi\,e^{-\xi t}\left[1-\int_0^t\widehat g_{j,r}(\tau)\,{\rm d}\tau\right],
 \qquad t>0,
 \label{equation:21}
\end{equation}
\begin{equation}
 \gamma_{j,r}(\lambda)
 ={\lambda\over \lambda+\xi}\,\widehat\gamma_{j,r}(\lambda+\xi)
 +{\xi\over \lambda+\xi},
 \qquad \lambda>0.
 \label{equation:25}
\end{equation}
\end{prop}
\begin{proof}
Making use of Eq.\ (\ref{equation:68}) and (\ref{equation:69}) for $k=r$, from 
(\ref{equation:19}) we obtain:
\begin{eqnarray}
 \int_0^{t} g_{j,r}(\tau)\,{\rm d}\tau \!\!\!\!
 &=& \!\!\!\! 
 1-e^{-\xi t}\sum_{n=r+1}^{+\infty}\widehat A_{j,n}^{\langle r\rangle}(t) 
 \nonumber \\
 &=& \!\!\!\! 
 1-e^{-\xi t}\left[1-\int_0^t \widehat g_{j,r}(\tau)\,{\rm d}\tau\right].
 \label{equation:23}
\end{eqnarray}
Differentiating both sides of (\ref{equation:23}) with respect to $t$ 
we are led to (\ref{equation:21}). Taking the Laplace transform of both sides 
of (\ref{equation:21}) and making use of Fubini's theorem we then have:
$$
 \gamma_{j,r}(\lambda)
 =\widehat\gamma_{j,r}(\lambda+\xi)+{\xi\over \lambda+\xi}
 -\xi\int_0^{\infty}\widehat g_{j,r}(\tau) 
 \left[\int_{\tau}^{\infty} e^{-(\lambda+\xi)t}\,{\rm d}t\right]{\rm d}\tau,
 \qquad \lambda>0,
$$
from which Eq.\ (\ref{equation:25}) follows. 
\qed
\end{proof}
\par
In the limit as $t\to +\infty$, from (\ref{equation:23}) we obtain 
$\int_0^{+\infty} g_{j,r}(\tau)\,{\rm d}\tau =1$, implying that the first 
visit of $N(t)$ in state $r$ occurs with probability 1, whereas for 
$\widehat N(t)$ such probability may be less than 1. 
\par
A characterization of the distribution of the first-visit time 
in state $r$ is provided by the following 
\begin{thm}\label{theorem:2}
Let $Z$ be an exponentially distributed random variable independent of 
$\widehat T_{j,r}$ with mean $\xi^{-1}$. If $\Prob(\widehat T_{j,r}<+\infty)=1$, 
then for $j\in\{r+1,r+2,\ldots\}$ the r.v. 
\begin{equation}
 \Theta_{j,r}:=\min\big\{\widehat T_{j,r}, Z\big\}
 \label{equation:24}
\end{equation}
has the same distribution as $T_{j,r}$. 
\end{thm}
\begin{proof}
From (\ref{equation:24}) we have: 
$$
 \Prob(\Theta_{j,r}\leq t)
 =1-e^{-\xi t}+e^{-\xi t}\,P\big(\widehat T_{j,r}\leq t\big) 
 =\Prob(T_{j,r}\leq t),
 \qquad t\geq 0,
$$
where the last equality follows from (\ref{equation:23}). 
This completes the proof. 
\qed
\end{proof}
%
\subsection{First visit to state $k$}
We shall investigate the first-visit-time problem to any state 
$k\in\Stati$. Note that the case $j=r$ and $k=r+1$ is trivial and 
$$
 g_{r,r+1}(t)=\alpha_r\,e^{-\alpha_r t}, 
 \qquad 
 A_{r,r}^{\langle r+1\rangle}(t)=e^{-\alpha_r t},
 \qquad  t>0.
$$
In the following theorem $\gamma_{j,k}(\lambda)$ is expressed 
in terms of $\widehat\gamma_{j,k}(\lambda)$. 
\begin{thm}
For all $\lambda>0$ and for $j,k\in\Stati$, $j\neq k$, there results: 
\begin{equation}
 \gamma_{j,k}(\lambda)
 =\frac{\lambda\,\widehat\gamma_{j,k}(\lambda+\xi)+\xi\,\widehat\gamma_{r,k}(\lambda+\xi)}
 {\lambda+\xi\,\widehat\gamma_{r,k}(\lambda+\xi)}.
 \label{equation:54}
\end{equation}
\end{thm}
\begin{proof}
For all $t>0$ and for $j,k\in\Stati$, $j\neq k$, the following 
renewal equation holds: 
\begin{equation}
 p_{j,k}(t)=\int_0^t g_{j,k}(\tau)\,p_{k,k}(t-\tau)\,{\rm d}\tau,
 \label{equation:28} 
\end{equation}
Taking the Laplace transform of (\ref{equation:28}) and recalling (\ref{equation:6}) 
we have:
$$
 \gamma_{j,k}(\lambda)
 =\frac{\pi_{j,k}(\lambda)}{\pi_{k,k}(\lambda)}
 =\frac{\lambda\,\widehat \pi_{j,k}(\lambda+\xi)
 +\xi\,\widehat \pi_{r,k}(\lambda+\xi)}
 {\lambda\,\widehat \pi_{k,k}(\lambda+\xi)
 +\xi\,\widehat \pi_{r,k}(\lambda+\xi)},
$$
that leads one to (\ref{equation:54}) by virtue of the relation 
$\widehat\gamma_{j,k}(\lambda)=\widehat\pi_{j,k}(\lambda)/\widehat\pi_{k,k}(\lambda)$, 
$\lambda>0$. 
\qed
\end{proof}
\par
Setting $\lambda=0$ in Eq.\ (\ref{equation:54}) yields $\gamma_{j,k}(0)=1$, so 
that the first visit of $N(t)$ to any state $k$ occurs with unit probability, 
whereas for $\widehat N(t)$ this event may occur with probability less than 1. 
We note that Eq.\ (\ref{equation:54}) reduces to (\ref{equation:25}) if $k=r$. 
\par
Making use of Eq.\ (\ref{equation:54}) we can now immediately express mean and 
variance of $T_{j,k}$ in terms of the Laplace transform of $\widehat g_{j,r}(t)$, 
for all $j,k\in\Stati$, $j\neq k$: 
\begin{eqnarray}
 && \EE(T_{j,k}) = \frac{1-\widehat\gamma_{j,k}(\xi)}{\xi\,\widehat\gamma_{r,k}(\xi)},
 \label{equation:31} \\
 && \Var(T_{j,r}) = \frac{1}{\xi^2\,\widehat\gamma_{r,k}^2(\xi)}
 \left\{1-\widehat\gamma_{j,k}^2(\xi)
 +2\xi\,\big[1-\widehat\gamma_{j,k}(\xi)\big]
 \frac{\rm d}{{\rm d}\xi}\widehat\gamma_{r,k}(\xi)\right.
 \nonumber \\
 && \hspace{2cm}
 \left. +2\xi\,\widehat\gamma_{r,k}(\xi)\frac{\rm d}{{\rm d}\xi}\widehat\gamma_{j,k}(\xi)\right\}.
 \label{equation:32}
\end{eqnarray}
\par
Eqs.\ (\ref{equation:31}) and (\ref{equation:32}) will be used in the Appendix 
for the analysis of an immigration-emigration process with catastrophes. 
\par
Next proposition yields the Laplace transform of the $k$-avoiding transition 
probability $A^{\langle k\rangle}_{j,n}(t)$. 
\begin{prop}\label{prop:1}
Under each of the following mutually exclusive assumptions: \\
(i) \ $k\in\{r+1,r+2,\ldots\}$ and $j,n\in\{r,r+1,\ldots,k-1\}$, \\
(ii) \ $n,k\in\Stati$, $n\neq k$, and $j\in\{k+1,k+2,\ldots\}$, \\
for all $\lambda>0$ we have: 
\begin{eqnarray}
 && \hspace{-0.5cm}
 {\cal A}^{\langle k\rangle}_{j,n}(\lambda)
 =\widehat\pi_{j,n}(\lambda+\xi)-\widehat\gamma_{j,k}(\lambda+\xi)\,\widehat\pi_{k,n}(\lambda+\xi)
 \nonumber \\
 && 
 \hspace{1cm}
 +\,\xi\;\frac{1-\widehat\gamma_{j,k}(\lambda+\xi)}
 {\lambda +\xi\,\widehat\gamma_{r,k}(\lambda+\xi)}
 \,\big[\widehat\pi_{r,n}(\lambda+\xi)
 -\widehat\gamma_{r,k}(\lambda+\xi)\,\widehat\pi_{k,n}(\lambda+\xi)\big].  
 \label{equation:52}
\end{eqnarray}
\end{prop}
\begin{proof}
Under assumptions (i) or (ii) there holds: 
$$
 A_{j,n}^{\langle r\rangle}(t)
 =p_{j,n}(t)-\int_{0}^t g_{j,k}(\tau)\,p_{k,n}(t-\tau)\,{\rm d}\tau, 
 \qquad t>0.
$$
Taking the Laplace transform of this equation and recalling 
(\ref{equation:6}) and (\ref{equation:54}) we obtain: 
\begin{eqnarray*}
 && {\cal A}^{\langle k\rangle}_{j,n}(\lambda)
 =\widehat\pi_{j,n}(\lambda+\xi)+\frac{\xi}{\lambda}\,\widehat\pi_{r,n}(\lambda+\xi) \\
 && \hspace{1.5cm}
 -\frac{\lambda\,\widehat\gamma_{j,k}(\lambda+\xi)+\xi\,\widehat\gamma_{r,k}(\lambda+\xi)}
 {\lambda+\xi\,\widehat\gamma_{r,k}(\lambda+\xi)}
 \,\left[\widehat\pi_{k,n}(\lambda+\xi)+\frac{\xi}{\lambda}\,\widehat\pi_{r,n}(\lambda+\xi)\right], 
\end{eqnarray*}
that leads us to (\ref{equation:52}) after some calculations.
\qed
\end{proof}
\par
We point out that under assumption (i) of Proposition \ref{prop:1}, from the 
well-known relation 
$$
 \widehat A_{j,n}^{\langle k\rangle}(t)
 =\widehat p_{j,n}(t)-\int_0^t \widehat g_{j,k}(\tau)\,\widehat p_{k,n}(t-\tau)\,{\rm d}\tau,
 \qquad t>0,
$$
we obtain:
$$
 \widehat{\cal A}^{\langle k\rangle}_{j,n}(\lambda)
 =\widehat\pi_{j,n}(\lambda)-\widehat\gamma_{j,k}(\lambda)\,\widehat\pi_{k,n}(\lambda),
 \qquad \lambda>0.
$$
Eq.\ (\ref{equation:52}) can thus be also expressed as 
$$
 {\cal A}^{\langle k\rangle}_{j,n}(\lambda)
 =\widehat{\cal A}^{\langle k\rangle}_{j,n}(\lambda+\xi) 
 +\xi\;\frac{1-\widehat\gamma_{j,k}(\lambda+\xi)}
 {\lambda +\xi\,\widehat\gamma_{r,k}(\lambda+\xi)}
 \, \widehat{\cal A}^{\langle k\rangle}_{r,n}(\lambda+\xi) ,  
 \qquad j,n\in\{r,r+1,\ldots,k-1\}.
$$
We remark that since $N(t)$ is not skip-free to the left, 
${\cal A}^{\langle k\rangle}_{j,n}(\lambda)$ cannot be expressed in terms 
of $\widehat {\cal A}^{\langle k\rangle}_{j,n}(\lambda)$ when $j>k$. 
Moreover, since $N(t)$ is skip-free to the right, for all $t>0$ 
and $j\in\{r,r+1,\ldots, k-1\}$ the following equation holds:
\begin{equation}
 g_{j,k}(t)=\int_0^t g_{j,n}(\tau)\,g_{n,k}(t-\tau)\,{\rm d}\tau,
 \qquad n=j+1,j+2,\ldots,k-1.
 \label{equation:30}
\end{equation}
Instead, Eq.\ (\ref{equation:30}) does not hold when $j>k$ because $N(t)$ 
is not skip-free to the left. 
\section{First occurrence of effective catastrophe}
This section focuses on ``effective'' catastrophes, i.e.\ on catastrophes 
that are able to change the state of $N(t)$. Therefore, the occurrence of 
catastrophes while the process is in state $r$ is not taken into account. 
Consequently, the catastrophes first-occurrence time is no longer 
exponentially distributed. 
\par
Let us denote by $C_{j,r}$ the first occurrence of an effective catastrophe, 
when $N(0)=j$, with $j\in\Stati$. An effective catastrophe, or shortly 
``a catastrophe'' from now on, produces a transition, with rate $\xi$, 
from any state $n>r$ to the reflecting state $r$. Hence, certain 
transitions from $r+1$ to $r$ may be due to the occurrence of a 
catastrophe (with rate $\xi$), whereas the remaining transitions 
are due to a death (with rate $\beta_{r+1}$). 
\par
In order to investigate on the features of $C_{j,r}$, let us refer to a modified 
birth-death process with catastrophes that will be denoted as $\{M(t);\, t\geq 0\}$. 
This is assumed to be defined on the state-space $\{r-1,r,r+1,\ldots\}$. Its behavior 
is identical to that of $N(t)$, the only difference being that the effect of a 
catastrophe from state $n>r$ is a jump from $n$ to the absorbing state $r-1$. 
In other words, the allowed transitions are the following: \\
(i) \ $n\to n+1$ with rate $\alpha_n$, for $n=r,r+1,\ldots$, \\
(ii) \ $n\to n-1$ with rate $\beta_n$, for $n=r+1,r+2,\ldots$, \\
(iii) \ $n\to r-1$ with rate $\xi$, for $n=r+1,r+2,\ldots$. 
\par
For all $t\geq 0$ and $j\in\Stati$, $n\in\{r-1,r,r+1,\ldots\}$, let us now 
consider the transition probabilities of the modified process 
$$
 h_{j,n}(t)=\Prob\{M(t)=n\,|\,M(0)=j\}.
$$
The link between $M(t)$ and $C_{j,r}$ is evident by noting that the transitions of 
$M(t)$ from $n>r$ to $r-1$ corresponds to the transitions of $N(t)$ 
from $n$ to $r$ due to a catastrophe. Hence, denoting by $d_{j,r}(t)$ the density 
of $C_{j,r}$ for all $t>0$ we have:
\begin{equation}
 \Prob(C_{j,r}<t)\equiv \int_0^t d_{j,r}(\tau)\,{\rm d}\tau = h_{j,r-1}(t), 
 \qquad j\in\Stati.
 \label{equation:38}
\end{equation}
Moreover, for all $j\in\Stati$ the following system of forward equations holds:
\begin{eqnarray}
 && \hspace{-1cm} \frac{\rm d}{{\rm d}t} h_{j,r-1}(t) 
 =\xi\,[1-h_{j,r-1}(t)-h_{j,r}(t)], 
 \nonumber \\
 && \hspace{-1cm} \frac{\rm d}{{\rm d}t} h_{j,r}(t) 
 =-\alpha_r\,h_{j,r}(t)+\beta_{r+1}\,h_{j,r+1}(t), 
 \nonumber \\
 && \hspace{-1cm} \frac{\rm d}{{\rm d}t} h_{j,n}(t) 
 =-(\alpha_n+\beta_n+\xi)\,h_{j,n}(t)+\alpha_{n-1}\,h_{j,n-1}(t) 
 +\beta_{n+1}\,h_{j,n+1}(t), 
  \label{equation:34} \\
 && \hspace{7.5cm}  n=r+1,r+2,\ldots,
 \nonumber
\end{eqnarray}
with initial condition 
$$
 h_{j,n}(0)=\delta_{j,n}.
$$
Let us denote by $\eta_{j,n}(\lambda)$ the Laplace transform of $h_{j,n}(t)$. In the following 
theorem we shall express $\eta_{j,n}(\lambda)$ in terms of $\widehat\pi_{j,n}(\lambda)$. 
\begin{thm}
For all $j\in\Stati$ and $\lambda>0$ we have:
\begin{eqnarray}
 && \hspace{-1.5cm}
 \eta_{j,r-1}(\lambda) 
 = \frac{\xi}{\lambda + \xi}\,\left[\frac{1}{\lambda}
 -\frac{\widehat\pi_{j,r}(\lambda+ \xi)}{1-\xi\,\widehat\pi_{r,r}(\lambda+ \xi)}\right],
 \label{equation:36} \\
 && \hspace{-1.5cm}
 \eta_{j,n}(\lambda) 
 = \widehat\pi_{j,n}(\lambda+ \xi)+\xi\,\widehat\pi_{r,n}(\lambda+ \xi)\,
 \frac{\widehat\pi_{j,r}(\lambda+ \xi)}{1-\xi\,\widehat\pi_{r,r}(\lambda+ \xi)},
 \qquad n=r,r+1,\ldots.
 \label{equation:37}
\end{eqnarray}
\end{thm}
\begin{proof}
We shall prove the theorem by considering separately the two cases (a) $j=r$ and (b) $j>r$. 
\par
(a) \ Let $j=r$. Taking the Laplace transform of second and third equation in 
(\ref{equation:34}), we obtain:
\begin{eqnarray}
 && \hspace{-0.5cm}
 (\lambda+\alpha_r)\,\eta_{r,r}(\lambda)-1=\beta_{r+1}\,\eta_{r,r+1}(\lambda),
 \nonumber \\
 && \hspace{-0.5cm}
 (\lambda+\alpha_n+\beta_n+\xi)\,\eta_{r,r}(\lambda)
 =\alpha_{n-1}\,\eta_{r,n-1}(\lambda)+\beta_{n+1}\,\eta_{r,n+1}(\lambda), 
 \qquad n\geq r.
 \label{equation:39}
\end{eqnarray}
Recalling that $\pi_{r,n}(\lambda)$ is the Laplace transform of $p_{r,n}(t)$,
we look for a solution of the form
\begin{equation}
 \eta_{r,n}(\lambda)=A(\lambda)\,\pi_{r,n}(\lambda),
 \qquad n=r+1,r+2,\ldots,
 \label{equation:40}
\end{equation}
where $A(\lambda)$ does not depend on $n$. Substituting Eq.\ (\ref{equation:40}) 
in (\ref{equation:39}), we obtain: 
\begin{eqnarray}
 && \hspace{-1cm}
 A(\lambda)\,(\lambda+\alpha_r)\,\pi_{r,r}(\lambda)-1
 =\beta_{r+1}\,A(\lambda)\,\pi_{r,r+1}(\lambda),
 \nonumber \\
 && \hspace{-1cm}
 (\lambda+\alpha_n+\beta_n+\xi)\,\pi_{r,r}(\lambda)
 =\alpha_{n-1}\,\pi_{r,n-1}(\lambda)+\beta_{n+1}\,\pi_{r,n+1}(\lambda), 
 \qquad n>r.
 \label{equation:41}
\end{eqnarray}
Moreover, taking the Laplace transform of (\ref{equation:2}) it follows:
\begin{eqnarray}
 && \hspace{-1cm}
 (\lambda+\alpha_r+\xi)\,\pi_{r,r}(\lambda)-1
 =\beta_{r+1}\,\pi_{r,r+1}(\lambda)+\frac{\xi}{\lambda},
 \nonumber \\
 && \hspace{-1cm}
 (\lambda+\alpha_n+\beta_n+\xi)\,\pi_{r,r}(\lambda)
 =\alpha_{n-1}\,\pi_{r,n-1}(\lambda)+\beta_{n+1}\,\pi_{r,n+1}(\lambda), 
 \qquad n>r.
 \label{equation:42}
\end{eqnarray}
Comparing (\ref{equation:42}) with (\ref{equation:41}) we obtain: 
$$
 A(\lambda)=\frac{\lambda}{\lambda+\xi-\lambda\,\xi\,\pi_{r,r}(\lambda)},
$$
that does not depend on $n$. By virtue of (\ref{equation:40}) and 
(\ref{equation:6}), this yields Eq.\ (\ref{equation:37}) for $j=r$. 
\par
(b) \ Let $j>r$. Taking the Laplace transform of (\ref{equation:34}), we have: 
\begin{eqnarray*}
 && \hspace{-0.5cm}
 (\lambda+\alpha_r)\,\eta_{j,r}(\lambda)
 =\beta_{r+1}\,\eta_{j,r+1}(\lambda), \\
 && \hspace{-0.5cm}
 (\lambda+\alpha_n+\beta_n+\xi)\,\eta_{j,n}(\lambda)
 =\alpha_{n-1}\,\eta_{j,n-1}(\lambda)+\beta_{n+1}\,\eta_{j,n+1}(\lambda), 
 \qquad n>r, \;\; n\neq j, \\
 && \hspace{-0.5cm}
 (\lambda+\alpha_j+\beta_j+\xi)\,\eta_{j,j}(\lambda)-1
 =\alpha_{j-1}\,\eta_{j,j-1}(\lambda)+\beta_{j+1}\,\eta_{j,j+1}(\lambda). 
\end{eqnarray*}
Similarly to case (a) above there results: 
\begin{equation}
 \eta_{j,n}(\lambda)
 =B(\lambda)\,\pi_{j,n}(\lambda)+C(\lambda)\,\pi_{r,n}(\lambda),
 \qquad n=r-1,r,\ldots,
 \label{equation:44}
\end{equation}
where 
$$
 B(\lambda)=1, \qquad 
 C(\lambda)=\frac{\xi\,[\lambda\,\pi_{j,r}(\lambda)-1]}
 {\lambda+\xi-\xi\,\lambda\,\pi_{r,r}(\lambda)}.
$$
Hence, making use of (\ref{equation:44}) and (\ref{equation:6}), some straightforward 
calculations led us to Eq.\ (\ref{equation:37}) for $j>r$. Finally, taking the Laplace 
transform of the first equation in (\ref{equation:34}) we obtain:
$$
 \eta_{j,r-1}(t)
 =\frac{\xi}{\lambda+\xi}\left[\frac{1}{\lambda}-\eta_{j,r}(t)\right].
$$
Making again use of (\ref{equation:6}) Eq.\ (\ref{equation:36}) then follows. 
\qed
\end{proof}
\par
Let us now denote by $\delta_{j,r}(\lambda)$ the Laplace transform of $d_{j,r}(t)$, 
$j\in\Stati$. Taking the Laplace transform of both sides of (\ref{equation:38}) 
and making use of (\ref{equation:36}), for $j\in\Stati$ we have:
\begin{equation}
 \delta_{j,r}(\lambda)=\lambda\,\eta_{j,r-1}(\lambda) 
 =\frac{\xi}{\lambda+\xi}
 -\frac{\lambda}{\lambda+\xi}\,
 \frac{\xi\,\widehat\pi_{j,r}(\lambda+\xi)}{1-\xi\,\widehat\pi_{r,r}(\lambda+\xi)}.
 \label{equation:10}
\end{equation}
\par
The following renewal equation for $\widehat N(t)$
$$
 \widehat p_{j,r}(t)= 
 \int_0^t \widehat g_{j,r}(\tau)\,\widehat p_{r,r}(t-\tau)\,{\rm d}\tau,
 \qquad t>0,
$$
holding for all $j>r$, yields:
$$
 \widehat \pi_{j,r}(\lambda)
 =\widehat \gamma_{j,r}(\lambda)\,\widehat \pi_{r,r}(\lambda),
 \qquad \lambda>0.
$$
Hence, for $j>r$, (\ref{equation:10}) can be re-written as 
$$
 \delta_{j,r}(\lambda)
 =\frac{\xi}{\lambda+\xi}+\frac{\lambda}{\lambda+\xi}\,\widehat\gamma_{j,r}(\lambda+\xi)
 -\frac{\lambda}{\lambda+\xi}\,
 \frac{\widehat\gamma_{j,r}(\lambda+\xi)}{1-\xi\,\widehat\pi_{r,r}(\lambda+\xi)},
 \qquad \lambda>0.
$$
Making use of (\ref{equation:25}) in the right-hand-side, we have an alternative 
expression for $\delta_{j,r}(\lambda)$, for all $\lambda>0$ and $j\in\{r+1,r+2,\ldots\}$:
\begin{equation}
 \delta_{j,r}(\lambda)
 =\gamma_{j,r}(\lambda)-\frac{\lambda}{\lambda+\xi}\,
 \frac{\widehat\gamma_{j,r}(\lambda+\xi)}{ 1-\xi\,\widehat \pi_{r,r}(\lambda+\xi)}.
 \label{equation:45}
\end{equation}
\begin{prop}
For all $j\in\Stati$ there holds: 
\begin{equation}
 \EE(C_{j,r}) = \frac{1}{\xi}\,
 \left[\frac{\widehat\gamma_{j,r}(\xi)}{1-\xi\,\widehat\pi_{r,r}(\xi)}\right]. 
 \label{equation:46} 
\end{equation}
\end{prop}
\begin{proof}
For $j>r$, (\ref{equation:46}) follows by differentiating the right-hand-side 
of (\ref{equation:45}) with respect to $\lambda$ and then setting $\lambda=0$. 
The same procedure starting from (\ref{equation:10}) yields (\ref{equation:46}) 
in the case $j=r$. 
\qed
\end{proof}
%
\section{An extension to time-non-homogeneous processes}
In order to extend some of the above results to the time-non-homogeneous case, in 
this section $\{N(t);\, t\geq t_0\}$, $t_0\geq 0$, will denote a time non-homogeneous 
birth-death process with catastrophes, defined on the state-space $\Stati=\{r,r+1,r+2,\ldots\}$, 
where births occur with rates $\alpha_n(t)$ and deaths with rates $\beta_n(t)$. Moreover, 
catastrophes are assumed to occur according to a non-homogeneous Poisson process 
characterized by intensity function $\xi(t)$. We assume that $\alpha_n(t)$, 
$\beta_n(t)$ and $\xi(t)$ are continuous bounded functions such that 
$\alpha_n(t)>0$, $\beta_n(t)\geq 0$ and $\xi(t)>0$ for all $t\geq t_0$, 
with $\int_{t_0}^{+\infty}\xi(t)\,{\rm d}t=+\infty$. 
Moreover, $\{\widehat N(t);\, t\geq t_0\}$ will denote the time non-homogeneous 
birth-death process obtained from $N(t)$ by removing the possibility of catastrophes, 
i.e.\ by setting $\xi(t)=0$ for all $t\geq t_0$. For all $t\geq t_0$ and $j,n\in\Stati$ 
the transition probabilities of $N(t)$ and $\widehat N(t)$ will be denoted by
\begin{equation}
 p_{j,n}(t\,|\,t_0)=\Prob\{N(t)=n\,|\,N(t_0)=j\}, 
 \qquad  
 \widehat p_{j,n}(t\,|\,t_0)
 =\Prob\{\widehat N(t)=n\,|\,\widehat N(t_0)=j\}. 
 \label{equation:4}
\end{equation}
Assuming that $\widehat N(t)$ is non-explosive, i.e.\ such that 
$\sum_{n=r}^{+\infty}\widehat p_{j,n}(t\,|\,t_0)=1$ for all $j\in\Stati$ and 
$t\geq t_0$, we shall extend Eqs.\ (\ref{equation:5}) and (\ref{equation:27}) to the 
time-non-homogeneous case. To this end, we note that by making use of the forward equations 
for probabilities (\ref{equation:4}), for all $j,n\in\Stati$ and $t>t_0$ there hold: 
\begin{equation}
 p_{j,n}(t\,|\,t_0)=\exp\Big\{-\int_{t_0}^t \xi(u)\,{\rm d}u\Big\}\,\widehat p_{j,n}(t\,|\,t_0)
 +\int_{t_0}^{t}\xi(\tau)\,\exp\Big\{-\int_{\tau}^t \xi(u)\,{\rm d}u\Big\}\,
 \widehat p_{r,n}(t\,|\,\tau)\,{\rm d}\tau,
 \label{equation:65}
\end{equation}
and
\begin{eqnarray*}
 && \hspace{-1cm} 
 \EE\big[N(t)\,|\,N(t_0)=j\big]
 =\exp\Big\{-\int_{t_0}^t \xi(u)\,{\rm d}u\Big\}\,
 \EE\big[\widehat N(t)\,|\,\widehat N(t_0)=j\big]   \\
 && \hspace{2.2cm}
 +\int_{t_0}^{t}\xi(\tau)\,\exp\Big\{-\int_{\tau}^t \xi(u)\,{\rm d}u\Big\}\,
 \EE\big[\widehat N(t)\,|\,\widehat N(\tau)=r\big]\,{\rm d}\tau. 
\end{eqnarray*}
Moreover, similarly to the cases of Eqs.\ (\ref{equation:19}) and (\ref{equation:21}),  
it is not hard to prove that the $r$-avoiding transition probabilities and the 
pdf's of the first-visit time to state $r$ in the time-non-homogeneous case 
are related as follows for all $t\geq t_0$ and $j\in\{r+1, r+2,\ldots\}$:
\begin{eqnarray}
 && A_{j,n}^{\langle r\rangle}(t\,|\,t_0)
 =\exp\Big\{-\int_{t_0}^t \xi(u)\,{\rm d}u\Big\}\,
 \widehat A_{j,n}^{\langle r\rangle}(t\,|\,t_0), 
 \qquad n\in\{r+1,r+2,\ldots \}, 
 \nonumber \\
 && g_{j,r}(t\,|\,t_0)
 =\exp\Big\{-\int_{t_0}^t\xi(u)\,{\rm d}u\Big\}\,\widehat g_{j,r}(t\,|\,t_0) 
 \nonumber \\
 && \hspace{1.6cm}
 +\,\xi(t)\,\exp\Big\{-\int_{t_0}^t\xi(u)\,{\rm d}u\Big\}
 \Big[1-\int_{t_0}^t\widehat g_{j,r}(\tau\,|\,t_0)\,{\rm d}\tau\Big]. 
 \label{equation:71}
\end{eqnarray}
\par
Let $\widehat T_{j,r}$ denote the first-visit time in state $r$ for the 
time-non-homogeneous birth-death process $\widehat N(t)$. The characterization 
given in Theorem \ref{theorem:2} can be extended to this time-non-homogeneous case. 
Indeed, let $Z$ be a random variable with hazard function $\xi(t)$, $t\geq t_0$, 
i.e.\ with distribution function 
$F(t)=1-\exp\{-\int_{t_0}^t\xi(u)\,{\rm d}u\}$, $t\geq t_0$. If $Z$ is 
independent of $\widehat T_{j,r}$, and if $\Prob(\widehat T_{j,r}<+\infty)=1$, 
then it is possible to prove that for $j\in\{r+1,r+2,\ldots\}$ the r.v.\ 
(\ref{equation:24}) has the same distribution as the first-visit time $T_{j,r}$. 
\section{Some additional results}
Let us denote by $N$ the random variable that describes the steady state 
of $N(t)$ in the time-homogeneous case, and let 
\begin{equation}
 q_n:=\Prob(N=n)=\lim_{t\to +\infty}p_{j,n}(t), 
 \qquad j,n\in\Stati.
 \label{equation:7}
\end{equation}
Then, from (\ref{equation:2}) we obtain:
\begin{equation}
 \begin{array}{l}
 -(\alpha_r+\xi)\,q_{r}+\beta_{r+1}\,q_{r+1}+\xi=0, \\
 -(\alpha_n+\beta_n+\xi)\,q_{n}+\alpha_{n-1}\,q_{n-1}+\beta_{n+1}\,q_{n+1}=0, 
 \qquad n=r+1,r+2,\ldots.
 \end{array}
 \label{equation:61}
\end{equation}
Making use of a Tauberian theorem, from Eq.\ (\ref{equation:6}) there results: 
\begin{equation}
 q_n=\xi\,\widehat \pi_{r,n}(\xi), 
 \qquad t\geq 0, \;\; n\in\Stati.
 \label{equation:8}
\end{equation}
\par
From the assumed non explosivity of $\widehat N(t)$ and from (\ref{equation:8}) it follows 
that $N(t)$ possesses a steady-state distribution, with $\sum_{n=r}^{+\infty}q_n=1$. 
We stress that due to (\ref{equation:8}), such a distribution can be obtained from 
the Laplace transform of the transition probabilities in absence of catastrophes. 
Moreover, we underline that (\ref{equation:8})'s satisfy system (\ref{equation:61}). 
\par
In Eq.\ (\ref{equation:5}) we have expressed the transition probability of $N(t)$ 
in terms of that of $\widehat N(t)$. In order to provide a probabilistic 
interpretation of $p_{r,n}(t)$, let us now introduce a time-homogeneous birth-death 
process $\{N^*(t);\, t\geq 0\}$ characterized by birth and death rates 
\begin{equation}
 \alpha_n^*=\alpha_n\,\displaystyle
 \frac{\displaystyle\sum_{k=n}^{+\infty}q_k}{\displaystyle\sum_{k=n+1}^{+\infty}q_k}  
 \qquad (n\geq r), 
 \qquad 
 \beta_n^*=\beta_n\,\displaystyle
 \frac{\displaystyle\sum_{k=n+1}^{+\infty}q_k}{\displaystyle\sum_{k=n}^{+\infty}q_k}   
 \qquad (n>r), 
 \label{equation:9}
\end{equation}
with state-space $\Stati$, where $r$ is a reflecting state. Assuming that 
$\Prob\{N^*(0)=r\}=1$, we denote by $p^*_{r,n}(t)$ the transition probabilities 
of $N^*(t)$. 
\begin{thm}\label{theorem:3}
Let $u_{r,n}(t)$ be the transition probabilities of 
\begin{equation}
 U(t):=\min\{N^*(t), N\}, 
 \qquad t\geq 0,
 \label{equation:11}
\end{equation}
where $N$ is independent of $N^*(t)$ and distributed as in (\ref{equation:7}). 
Then, 
\begin{equation}
 u_{r,n}(t)= p_{r,n}(t)  
 \label{equation:12}
\end{equation}
for all $t\geq 0$ and $n\in\Stati$.
\end{thm}
\begin{proof}
From Eq.\ (\ref{equation:11}) it follows 
\begin{equation}
 u_{r,n}(t)= q_{n}  \sum_{k=n}^{+\infty}p^*_{r,k}(t)
 + p^*_{r,n}(t)\sum_{k=n+1}^{+\infty}q_{k}, 
 \qquad n\in\Stati. 
 \label{equation:13}
\end{equation}
Differentiating both sides of (\ref{equation:13}) and making 
use of the forward equations of $p^*_{r,n}(t)$, we have: 
\begin{eqnarray}
 && \hspace{-1cm} \frac{\rm d}{{\rm d}t} u_{r,r}(t) 
 =-\alpha_r^*\,(1-q_r)\,p_{r,r}^*(t)+\beta_{r+1}^*\,(1-q_r)\,p_{r,r+1}^*(t), 
 \nonumber \\
 && \hspace{-1cm} \frac{\rm d}{{\rm d}t} u_{r,n}(t) 
 =-\left(\alpha_n^*\sum_{k=n+1}^{+\infty}q_{k}+\beta_n^*\sum_{k=n}^{+\infty}q_{k}\right)\,p_{r,n}^*(t) 
 +\alpha_{n-1}^*\,p_{r,n-1}^*(t) \sum_{k=n}^{+\infty}q_{k}
 \nonumber \\
 && \hspace{1cm} 
 +\,\beta_{n+1}^*\,p_{r,n+1}^*(t) \sum_{k=n+1}^{+\infty}q_{k},
 \hspace{1.5cm} n=r+1,r+2,\ldots.
  \label{equation:14}
\end{eqnarray}
Moreover, from (\ref{equation:13}) we obtain:   
\begin{eqnarray}
 && \hspace{-1cm} -(\alpha_r+\xi)\,u_{r,r}(t)+\beta_{r+1}\,u_{r,r+1}(t)+\xi 
 =-\alpha_r^*\,p^*_{r,r}(t)+\beta_{r+1}^*(1-q_r-q_{r+1})\,p^*_{r,r+1}(t), 
 \nonumber \\
 && \hspace{-1cm} -(\alpha_n+\beta_n+\xi)\,u_{r,n}(t)+\alpha_{n-1}\,u_{r,n-1}(t) 
 +\beta_{n+1}\,u_{r,n+1}(t)
 \nonumber \\ 
 && \hspace{0cm}  
 =\left[-(\alpha_n+\beta_n+\xi)\sum_{k=n}^{+\infty}q_{k}+\alpha_{n-1}q_{n-1}\right]p_{r,n}^*(t) 
 +\alpha_{n-1}\,p_{r,n-1}^*(t)\sum_{k=n-1}^{+\infty}q_{k}
 \nonumber \\ 
 && \hspace{0cm}  
 +\beta_{n+1}\,p_{r,n+1}^*(t)  \sum_{k=n+2}^{+\infty}q_{k},
 \hspace{3cm} n=r+1,r+2,\ldots.
  \label{equation:15}
\end{eqnarray}
Due to (\ref{equation:61}) and (\ref{equation:9}), the right-hand-sides of 
Eqs.\ (\ref{equation:14}) and (\ref{equation:15}) are identical, so that $u_{r,n}(t)$ 
satisfy the system of forward equations (\ref{equation:2}). Furthermore, due to 
(\ref{equation:11}) or (\ref{equation:13}), initial condition $u_{r,n}(0)=\delta_{r,n}$ 
holds. Hence, (\ref{equation:12}) follows. 
\qed
\end{proof}
\par
An immediate consequence of Theorem \ref{theorem:3} is that $p_{r,n}(t)$ can be expressed 
in terms of probabilities $p^*_{r,k}(t)$ and $q_k$. Indeed, due to (\ref{equation:12}) 
and (\ref{equation:13}), for all $t\geq 0$ and $n\in\Stati$ we have: 
\begin{equation}
 p_{r,n}(t)= q_{n}  \sum_{k=n}^{+\infty}p^*_{r,k}(t)
 + p^*_{r,n}(t)\sum_{k=n+1}^{+\infty}q_{k}. 
 \label{equation:16}
\end{equation}
\par
Hereafter we show that this result can be extended to the more general 
case of arbitrary initial state. 
\begin{thm}\label{theorem:4}
For all $t\geq 0$ and $j,n\in\Stati$ we have:
\begin{equation}
 p_{j,n}(t)
 = e^{-\xi t}\left[\widehat p_{j,n}(t)-\widehat p_{r,n}(t)\right]
 +q_{n}  \sum_{k=n}^{+\infty}p^*_{r,k}(t)
 + p^*_{r,n}(t)\sum_{k=n+1}^{+\infty}q_{k}. 
 \label{equation:17}
\end{equation}
\end{thm}
\begin{proof}
Comparing Eq.\ (\ref{equation:5}) written for $j=r$ with (\ref{equation:16}) we obtain:  
\begin{eqnarray*}
 \xi\int_0^{t}e^{-\xi \tau}\,\widehat p_{r,n}(\tau)\,{\rm d}\tau \!\!\!\!
 &=& \!\!\!\! p_{r,n}(t)-e^{-\xi t}\,\widehat p_{r,n}(t) \\
 &=& \!\!\!\! q_{n}  \sum_{k=n}^{+\infty}p^*_{r,k}(t)
 + p^*_{r,n}(t)\sum_{k=n+1}^{+\infty}q_{k}
 -e^{-\xi t}\,\widehat p_{r,n}(t).
\end{eqnarray*}
Making use of this equation in (\ref{equation:5}) we are finally 
led to (\ref{equation:17}). 
\qed
\end{proof}
%
\section*{Appendix}
Hereafter we shall apply the results obtained above to some birth-death 
processes with catastrophes of interest in biological contexts.
\subsection*{A1. Birth process with catastrophes}
Let $N(t)$ denote the number of individuals present at time $t$ in a birth process 
with catastrophes on $\Stati=\{r,r+1,\ldots\}$, with 
$$
 \alpha_n>0 \quad (n=r,r+1,\ldots),
 \qquad 
 \beta_n=0 \quad (n=r+1,r+2,\ldots). 
$$
We assume that 
$$
 \sum_{n=r}^{+\infty}\frac{1}{\alpha_n}=+\infty,
$$
so that 
$\sum_{n=r}^{+\infty}\widehat p_{r,n}(t)=1$ for all $t>0$ (cf.\ for instance 
Cox and Miller \cite{CoMi65}). Since
$$
 \widehat \pi_{r,n}(\lambda)
 =\left\{\begin{array}{ll}
 \ds{\frac{1}{\lambda+\alpha_r}}, & \quad n=r, \\
 \hfill & \hfill \\
 \ds{\frac{\alpha_r\alpha_{r+1}\cdots\alpha_{n-1}}
 {(\lambda+\alpha_r)(\lambda+\alpha_{r+1})\cdots(\lambda+\alpha_n)}}, & \quad n=r+1,r+2,\ldots, 
 \end{array}
 \right.
$$
we have: 
$$
 q_n=\xi\,\widehat \pi_{r,n}(\xi)
 =\left\{\begin{array}{ll}
 \ds{\frac{\xi}{\xi+\alpha_r}}, & \quad n=r, \\
 \hfill & \hfill \\
 \ds{\frac{\xi\alpha_r\alpha_{r+1}\cdots\alpha_{n-1}}
 {(\xi+\alpha_r)(\xi+\alpha_{r+1})\cdots(\xi+\alpha_n)}}, & \quad n=r+1,r+2,\ldots, 
 \end{array}
 \right.
$$
Note that \\
(i) \ if $\alpha_n=\alpha$, then the stationary distribution is geometric:
$$
 q_n=\frac{\xi}{\xi+\alpha}\,\left(\frac{\alpha}{\xi+\alpha}\right)^{n-r}, 
 \qquad n=r,r+1,\ldots,
$$
(ii) \ if $\alpha_n=\xi\,(n+k)$, with $k$ a positive integer, 
then the stationary distribution is  
$$
 q_n=\frac{r+k}{(n+k)(n+k+1)}, 
 \qquad n=r,r+1,\ldots.
$$
From (\ref{equation:25}) and (\ref{equation:10}) it follows   
$g_{j,r}(t)=d_{j,r}(t)=\xi\,e^{-\xi t}$, $t>0$. 
\subsection*{A2. Time non-homogeneous immigration-emigration process with catastrophes}
Let $N(t)$ be the continuous-time Markov chains with 
state-space $\{0,1,\ldots\}$ that describes the number of individuals in a 
population subject to an immigration-emigration process in the presence of 
catastrophes, with immigration rates $\alpha_n(t)=\alpha\,w(t)$, 
emigration rates $\beta_n(t)=\beta\,w(t)$ and catastrophe rate $\xi(t)$, 
with $\alpha>0$, $\beta>0$ and where $w(t)$ is a bounded, continuous and 
positive function such that $\int_{t_0}^{+\infty}w(t)\,{\rm d}t=+\infty$. 
Making use of some well-known results (see, for instance, Medhi \cite{Me91} and 
Bailey \cite{Ba57}), for $\widehat N(t)$ in the absence of catastrophes we have: 
\begin{eqnarray}
 && \hspace{-0.5cm} 
 \widehat p_{j,n}(t\,|\,t_0)
 =\exp\Big\{-(\alpha+\beta)\int_{t_0}^t w(u)\,{\rm d}u\Big\} 
 \nonumber \\
 && \hspace{0.4cm}
 \times \bigg[\rho^{(n-j)/2}\,I_{n-j}\Big(2\sqrt{\alpha\beta}\int_{t_0}^t w(u)\,{\rm d}u\Big) 
 +\rho^{(n-j-1)/2}\,I_{n+j+1}\Big(2\sqrt{\alpha\beta}\int_{t_0}^t w(u)\,{\rm d}u\Big)  
 \nonumber \\
 && \hspace{0.4cm}
 +(1-\rho)\rho^n\,\sum_{k=n+j+2}^{+\infty}\rho^{-k/2}\,
 I_k\Big(2\sqrt{\alpha\beta}\int_{t_0}^t w(u)\,{\rm d}u\Big)\bigg]
 \label{equation:26}
\end{eqnarray}
and: 
$$
 \widehat g_{j,0}(t\,|\,t_0)
 =\frac{j\,w(t)}{\int_{t_0}^t w(u)\,{\rm d}u}\,
 \exp\Big\{-(\alpha+\beta)\int_{t_0}^t w(u)\,{\rm d}u\Big\}\,\rho^{-j/2}\, 
 I_j\Big(2\sqrt{\alpha\beta}\int_{t_0}^t w(u)\,{\rm d}u\Big),
$$
where $\rho=\alpha/\beta$. Hence, the transition probabilities and the first-visit-time 
density of $N(t)$ can be easily obtained by making use of (\ref{equation:65}) and 
(\ref{equation:71}).  
\subsection*{A3. Immigration-emigration process with catastrophes}
Let $N(t)$ denote the number of individuals in a population subject to an 
immi\-gration-emigration process with catastrophes and characterized by constant 
birth and death rates $\alpha_n=\alpha$ and $\beta_n=\beta$, with state-space 
$\Stati=\{0,1,\ldots\}$. Since (see, for instance, Section~3.1 of Conolly~\cite{Co75}) 
$$
 \widehat \pi_{0,n}(\xi)
 =\frac{1}{\xi}\,(1-q)\,q^n, 
 \qquad 
 n=0,1,2,\ldots,
$$
with 
\begin{equation}
 q=\frac{\alpha+\beta+\xi-\sqrt{(\alpha+\beta+\xi)^2-4\alpha\beta}}{2\beta}, 
 \label{equation:49}
\end{equation}
Eq.\ (\ref{equation:8}) shows that the steady-state probabilities for 
the corresponding process $N(t)$ are given by:  
\begin{equation}
 q_n=(1-q)\,q^n,  \qquad 
 n=0,1,2,\ldots.
 \label{equation:62}
\end{equation}
By virtue of (\ref{equation:9}) and (\ref{equation:62}), the birth-death 
process $N^*(t)$ is characterized by the following rates: 
$$
 \alpha^*_n=\frac{\alpha}{q} \qquad (n\geq 0),
 \qquad
 \beta^*_n=\beta\,q \qquad (n\geq 1).
$$
Making use of (\ref{equation:17}) for all $t\geq 0$ and $j,n\in\{0,1,\ldots\}$, 
we have: 
\begin{equation}
 p_{j,n}(t)
 = e^{-\xi t}\left[\widehat p_{j,n}(t)-\widehat p_{0,n}(t)\right]
 +(1-q)q^{n}  \sum_{k=n}^{+\infty}p^*_{0,k}(t)
 + p^*_{0,n}(t)\,q^{n+1}, 
 \label{equation:63}
\end{equation}
where (see for instance Eq.\ (9.13) of Medhi \cite{Me91}, p.\ 120)
\begin{eqnarray*}
 && \hspace{-0.6cm} 
 p_{0,n}^*(t)=\exp\left\{-\left(\frac{\alpha}{q} +\beta\,q\right)t\right\} 
 \bigg[\left(\frac{\rho}{q^2}\right)^{n/2}\,I_n(2\sqrt{\alpha\beta}\, t) 
 +\left(\frac{\rho}{q^2}\right)^{(n-1)/2}\,I_{n+1}(2\sqrt{\alpha\beta}\, t)\\
 && \hspace{1cm}
 +\left(1-\frac{\rho}{q^2}\right)\left(\frac{\rho}{q^2}\right)^n\,
 \sum_{k=n+2}^{+\infty}\left(\frac{\rho}{q^2}\right)^{-k/2}\,I_k(2\sqrt{\alpha\beta}\, t)\bigg],
\end{eqnarray*}
with $\rho=\alpha/\beta$, and where $\widehat p_{j,n}(t)$ identifies with 
expression (\ref{equation:26}) upon setting $w(t)=1$. We note that Eq.\ 
(\ref{equation:63}) extends the results obtained in Section 2 of 
\cite{DiGiNoRi03} in which the case of initial state $j=0$ is treated. 
\par
Let us discuss the first-passage time of $N(t)$ in state $r$. 
Since $r=0$, in this case we have (see, for instance, Conolly~\cite{Co75}):  
\begin{equation}
 \widehat\gamma_{j,0}(\xi)
 =\left(\frac{\beta}{\alpha}\,q\right)^j
 =\left(\frac{\xi+\alpha+\beta-\sqrt{(\xi+\alpha+\beta)^2-4\alpha\beta}}{2\alpha}\right)^j,
 \label{equation:53}
\end{equation}
where $q$ is defined in (\ref{equation:49}). 
Making use of Eqs.\ (\ref{equation:31}) and (\ref{equation:32}) we thus obtain mean 
and variance of the first-visit time in state $0$ in the presence of catastrophes:
\begin{eqnarray}
 && \EE(T_{j,0})
 =\frac{1}{\xi}\left[1-\left(\frac{\beta}{\alpha}\,q\right)^j\right], \\
 && \Var(T_{j,0}) 
 =\frac{1}{\xi^2}\left[1-\left(\frac{\beta}{\alpha}\,q\right)^2
 +\frac{1}{\alpha}\;\frac{2\xi \beta q}{2\beta q-(\xi+\alpha+\beta)}\right].
\end{eqnarray}
Finally, recalling (\ref{equation:53}), from (\ref{equation:46}) we obtain 
the mean time of the first occurrence of an effective catastrophe: 
$$
 \EE(C_{j,0})
 =\frac{1}{\xi}+\frac{1-q}{q}\left(\frac{\beta}{\alpha}\,q\right)^j.
$$
\subsection*{A4. Immigration-death process with catastrophes}
Let $N(t)$ denote the number of individuals present at time $t$ in an immigrations-death 
process with catastrophes. Its state-space is $\{0,1,\ldots\}$, and its birth and 
death rates are  
$$
 \alpha_n=\nu \quad (n=0,1,\ldots),
 \qquad 
 \beta_n=\beta\,n \quad (n=1,2,\ldots), 
$$
with $\nu>0$ and $\beta>0$. Let us set $\rho=\nu/\beta$. 
It is well-known (cf.\ Cox and Miller \cite{CoMi65}) that 
$$
 \widehat p_{0,n}(t)
 =\frac{\left[\rho\left(1-e^{-\beta t}\right)\right]^n}{n!}\,
 \exp\left\{-\rho \left(1-e^{-\beta t}\right)\right\}, 
 \qquad n=0,1,\ldots\;.
$$
Then (see Eq.\ n.\ 3.383 of Gradshteyn and Ryzhik~\cite{GrRy80}, pag.\ 318),   
$$
 \widehat \pi_{0,n}(\lambda)
 =\frac{\rho^n}{n!}\,e^{-\rho}\,\frac{1}{\beta}\,B\left(n+1,\frac{\lambda}{\beta}\right)\, 
 \Phi\left(\frac{\lambda}{\beta},\frac{\lambda}{\beta}+n+1;\rho\right), 
 \qquad n=0,1,\ldots\;,
$$
where $B(a,b)=\Gamma(a)\,\Gamma(b)/\Gamma(a+b)$ is the Beta function and 
$$
 \Phi(a,c;x)=1+\sum_{n=1}^{+\infty}\frac{(a)_n}{(c)_n}\,\frac{x^n}{n!} 
$$
is the confluent hypergeometric function of first kind (Kummer function), 
with $(a)_n=a(a+1)\cdots(a+n-1)$. Hence, from (\ref{equation:8}) we have: 
\begin{equation}
 q_n=\xi\,\widehat \pi_{0,n}(\xi)
 =\frac{\rho^n}{n!}\,e^{-\rho}\,\frac{\xi}{\beta}\,B\left(n+1,\frac{\xi}{\beta}\right)\, 
 \Phi\left(\frac{\xi}{\beta}, \frac{\xi}{\beta}+n+1;\rho\right), 
 \qquad n=0,1,\ldots\;.
 \label{equation:55}
\end{equation}
In particular, the probability of having asymptotically zero individuals is: 
$$
 q_0=e^{-\rho}\,\frac{\xi}{\beta}\sum_{k=0}^{\infty}\frac{\rho^k}{k!}\,\frac{1}{k+\xi/\beta}.
$$
In the special case $\beta=\xi$, from (\ref{equation:55}) we have: 
$$
 q_n=\frac{1}{\rho}\,e^{-\rho}\sum_{i=n+1}^{+\infty}\frac{\rho^i}{i!},
 \qquad n=0,1,\ldots\;. 
$$
\subsection*{A5. Immigration-birth-death process with catastrophes}
Let $N(t)$ be the number of individuals present at time $t$ in an 
immigration-birth-death process with catastrophes, characterized by state-space 
$\{0,1,\ldots\}$ and by birth and death rates
$$
 \alpha_n=\alpha\,n+\nu \quad (n=0,1,\ldots),
 \qquad 
 \beta_n=\beta\,n \quad (n=1,2,\ldots), 
$$
with $\alpha>0$, $\nu>0$ and $\beta>0$. As is well-known, in absence of 
catastrophes for $n=0,1,\ldots$ there results: 
$$
 \widehat p_{0,n}(t)
 =\left\{
 \begin{array}{ll}
 \displaystyle\left(\frac{1}{1+\alpha t}\right)^{\frac{\nu}{\alpha}}\,
 \frac{1}{n!}\,\left(\frac{\nu}{\alpha}\right)_n\,
 \left(\frac{\alpha t}{1+\alpha t}\right)^n, & \alpha=\beta \\
 \hfill & \hfill \\
 \displaystyle\left[\frac{\alpha-\beta}
 {\alpha\,e^{(\alpha-\beta)t}-\beta}\right]^{\frac{\nu}{\alpha}}\,
 \frac{1}{n!}\,\left(\frac{\nu}{\alpha}\right)_n\,
 \left\{\frac{\alpha\left[e^{(\alpha-\beta)t}-1\right]}
 {\alpha\,e^{(\alpha-\beta)t}-\beta}\right\}^n, 
 & \alpha\neq\beta,
 \end{array}
 \right. 
$$
with 
$$
 \EE\big[\widehat N(t)\,|\,\widehat N(0)=j\big]
 =\left\{
 \begin{array}{ll}
 j+\nu\,t, & \alpha=\beta \\
 \hfill & \hfill \\
 j\,e^{(\alpha-\beta)t}+\ds{\frac{\nu\left[e^{(\alpha-\beta)t}-1\right]}{\alpha-\beta}}, 
 & \alpha\neq \beta.
 \end{array}
 \right.
$$
Let us now consider the process in the presence of catastrophes, whose transient 
probabilities have been studied in \cite{Sw97} by making use of a generating 
function approach. From (\ref{equation:27}) it follows:
\begin{equation}
 \EE\big[N(t)\,|\,N(0)=j\big]
 =\left\{
 \begin{array}{ll}
 j+\nu\,t, & \alpha=\beta+\xi \\
 \hfill & \hfill \\
 j\,e^{(\alpha-\beta-\xi)t}
 +\ds{\frac{\nu\left[e^{(\alpha-\beta-\xi)t}-1\right]}{\alpha-\beta-\xi}}, 
 & \alpha\neq \beta+\xi.
 \end{array}
 \right.
 \label{equation:56}
\end{equation}
It is interesting to note that the mean of $N(t)$ identifies with that of $\widehat N(t)$ 
if in Eq.\ (\ref{equation:56}) we substitute $\beta+\xi$ with $\beta$. Moreover, 
from (\ref{equation:56}) we have
$$
 \lim_{t\to+\infty}\EE\big[N(t)\,|\,N(0)=j\big]
 =\left\{
 \begin{array}{ll}
 \ds{\frac{\nu}{\beta+\xi-\alpha}}, & \alpha<\beta+\xi \\
 \hfill & \hfill \\
 +\infty, & \alpha\geq\beta+\xi.
 \end{array}
 \right.
$$
We shall now obtain the steady-state distribution of $N(t)$. We recall that in 
absence of catastrophes the steady-state distribution exists only if $\alpha<\beta$. 
\par
(i) Let $\alpha<\beta$. Then,  
%
\begin{eqnarray*}
 && \widehat \pi_{0,n}(\lambda)
 =\frac{1}{\alpha}\,\left(\frac{\alpha}{\beta}\right)^{n+1}\,
 \frac{1}{n!}\, \left(\frac{\nu}{\alpha}\right)_n\,
 \left(1-\frac{\alpha}{\beta}\right)^{\frac{\nu}{\alpha}-1}\,
 B\left(\frac{\lambda}{\beta-\alpha},n+1\right)\, \\
 && \hspace{1.5cm} 
 \times \,F\left(\frac{\nu}{\alpha}+n,\frac{\lambda}{\beta-\alpha};
 n+1+\frac{\lambda}{\beta-\alpha}; \frac{\alpha}{\beta}\right),
 \qquad n=0,1,\ldots,
\end{eqnarray*}
where: 
$$
 F(a,b;c;x)=1+\sum_{n=1}^{+\infty}\frac{(a)_n\,(b)_n}{(c)_n}\,\frac{x^n}{n!}.
$$
Hence, from (\ref{equation:8}) it follows 
\begin{eqnarray*}
 && q_n=\xi\,\widehat \pi_{0,n}(\xi)
 =\frac{\xi}{\alpha}\,\left(\frac{\alpha}{\beta}\right)^{n+1}\,
 \frac{1}{n!}\, \left(\frac{\nu}{\alpha}\right)_n\,
 \left(1-\frac{\alpha}{\beta}\right)^{\frac{\nu}{\alpha}-1}\,
 B\left(\frac{\xi}{\beta-\alpha},n+1\right)\, \\
 && \hspace{2.5cm} 
 \times \,F\left(\frac{\nu}{\alpha}+n,\frac{\xi}{\beta-\alpha};
 n+1+\frac{\xi}{\beta-\alpha}; \frac{\alpha}{\beta}\right), 
 \qquad n=0,1,\ldots\;.
\end{eqnarray*}
\par
(ii) If $\alpha=\beta$, then 
$$
 \widehat \pi_{0,n}(\lambda)
 =\frac{1}{\alpha}\,\left(\frac{\nu}{\alpha}\right)_n\,
 \left(\frac{\lambda}{\alpha}\right)^{\frac{\nu}{\alpha}-1}
 \,\Psi\left(\frac{\nu}{\alpha}+n,\frac{\nu}{\alpha};\frac{\lambda}{\alpha}\right),
$$
where 
$$
 \Psi(a,c;x)=\frac{1}{\Gamma(a)}\,\int_0^{+\infty}e^{-x t}\, t^{a-1}\,(1+t)^{c-a-1}\,{\rm d}t,
 \qquad \Re e(a)>0
$$ 
is the confluent hypergeometric function of the second kind. 
Since $\Psi(a-c+1,2-c;x)=x^{c-1}\,\Psi(a,c;x)$, 
from (\ref{equation:8}) we thus have: 
$$
 q_n=\xi\,\widehat \pi_{0,n}(\xi)
 = \left(\frac{\xi}{\alpha}\right)^{\frac{\nu}{\alpha}}\,\left(\frac{\nu}{\alpha}\right)_n\,
 \Psi\left(\frac{\nu}{\alpha}+n,\frac{\nu}{\alpha};\frac{\xi}{\alpha}\right), 
 \qquad n=0,1,\ldots\;.
$$
In particular, from identity $\Psi(a,a;x)=e^{x}\,\Gamma(1-a,x)$, 
we obtain: 
$$
 q_0=\left(\frac{\xi}{\alpha}\right)^{\frac{\nu}{\alpha}}\,e^{\xi/\alpha}\,
 \Gamma\left(1-\frac{\nu}{\alpha},\frac{\xi}{\alpha}\right).
$$
\par
(iii) Finally, if $\alpha>\beta$ there results: 
\begin{eqnarray*}
 && \widehat \pi_{0,n}(\lambda)
 =\frac{1}{\alpha}\,\frac{1}{n!}\, \left(\frac{\nu}{\alpha}\right)_n\,
 \left(1-\frac{\beta}{\alpha}\right)^{\frac{\nu}{\alpha}-1}\,
 B\left(\frac{\lambda}{\alpha-\beta}+\frac{\nu}{\alpha},n+1\right)\, \\
 && \hspace{1.5cm} 
 \times \,F\left(\frac{\nu}{\alpha}+n,\frac{\lambda}{\alpha-\beta}+\frac{\nu}{\alpha};
 \frac{\lambda}{\alpha-\beta}+\frac{\nu}{\alpha}+n+1; \frac{\beta}{\alpha}\right),
 \qquad n=0,1,\ldots,
\end{eqnarray*}
that by virtue of (\ref{equation:8}) implies:
\begin{eqnarray*}
 && \hspace{-0.5cm}  q_n=\xi\,\widehat \pi_{0,n}(\xi)
 =\frac{\xi}{\alpha}\,\frac{1}{n!}\, \left(\frac{\nu}{\alpha}\right)_n\,
 \left(1-\frac{\beta}{\alpha}\right)^{\frac{\nu}{\alpha}-1}\,
 B\left(\frac{\xi}{\alpha-\beta}+\frac{\nu}{\alpha},n+1\right)\, \\
 && \hspace{1.8cm} 
 \times \,F\left(\frac{\nu}{\alpha}+n,\frac{\xi}{\alpha-\beta}+\frac{\nu}{\alpha};
 \frac{\xi}{\alpha-\beta}+\frac{\nu}{\alpha}+n+1; \frac{\beta}{\alpha}\right), 
 \qquad n=0,1,\ldots .
\end{eqnarray*}
\subsection*{\bf Acknowledgments}
Work performed within a joint cooperation agreement between
Japan Science and Technology Corporation (JST) and Universit\`a di Napoli
Federico II, under partial support by INdAM (G.N.C.S). 

%
\end{document}